\documentclass[12pt]{amsart}

\newcommand{\eps}{\varepsilon}
\newcommand{\pa}{\partial}

\usepackage{amsmath,amsthm,amscd}
\newcommand{\frk}[1]{{\mathfrak{#1}}}

\title
[]{A Note on the Analytic Families of Compact Submanifolds of Complex
Manifolds}
\author[]{
Zhiqin Lu}
\date{January 10, 1998 revised version \today}
\address[Zhiqin Lu]
{Department of Mathematics\\
Columbia University\\
NY, NY 10027}
\email{zhiqin@math.columbia.edu}

\newtheorem{theorem}{Theorem}[section]
\newtheorem{lemma}{Lemma}[section]
\newtheorem{cor}{Corollary}[section]
\newtheorem{prop}{Proposition}[section]

\newtheorem{definition}{Definition}[section]

\theoremstyle{remark}

\begin{document}
\maketitle

\numberwithin{equation}{section}

\section{Introductions}

It is a classical theorem of Kodaira~\cite{K1} that if the first
cohomological group
of the normal bundle is zero, then the deformation of a
compact complex  submanifold
within the ambient complex manifold is unobstructed. However, it is
usually difficult
to check if such a cohomological group is indeed zero. Furthermore, 
as showed in \S 2, in
some cases, it will never be zero.

In this short note, we are going to prove: if the deformation of the
complex structure 
of a  compact complex manifold $M$
is unobstructed in the sense of Kodaira and Spencer,
and if $M\rightarrow V$ is an embedding 
to the complex manifold $V$
and $H^1(M,T_V|_M)=0$,
then any fiber in a neighborhood of the universal deformation space
$U$ at $M$ can be embedded holomorphically
to $V$.

Contrary to the case of  normal bundle,
it will be relatively easy to check the vanishing of the group 
$H^1(M,T_V|_M)$. For example, if the curvature of the manifold $V$ has
some kinds of  positivity along $M$, then the group vanishes.

The typical examples of $M$ are compact Calabi-Yau manifolds. Those
manifolds admit K\"ahler metric with zero Ricci curvature. 
In Tian~\cite{T1}, it is  proved that the deformation 
of the complex structure of a Calabi-Yau
manifold is unobstructed. The moduli space of a polarized 
Calabi-Yau manifold
is then a complex orbifold.

We use the similar method as that of Kodaira~\cite{K1}.
 That is, we construct a formal power series which gives the map we want.
Then we prove the  convergence of the power series. 
 Fortunately, in
our
case, the original method in ~\cite{K1} of proving the convergence can
also be applied
here.

 {\bf Acknowledgment.}
The author thanks his advisor, Professor G. Tian for the
help during the preparation of this paper.
He also thanks the referee for pointing out an error at
the early version of this paper.

\section{Analytic families of compact submanifolds}
Suppose $M$ is a compact complex manifold. We assume that the deformation
of
the complex structure of $M$ is unobstructed. That is, the 
universal deformation space $U$ of
$M$ is a complex manifold near $M$ with 
complex dimension $\dim\, H^1(M, T_M)$,
where $T_M$ denotes the holomorphic tangent bundle of $M$.

If $M$ is a Calabi-Yau manifold, then the deformation of the complex
structures on $M$ is unobstructed by a theorem of Tian~\cite{T1}.

We use the notations and definitions in ~\cite{K1}.

\begin{definition}
Suppose $N$ is a
complex manifold of dimension $r+n$.
By an analytic family of compact submanifolds of dimension $n$ of $N$ we
shall
mean a pair $({\mathcal M},U)$ of a complex manifold $U$ and a
complex analytic submanifold ${\mathcal M}$ of $N\times U$ of
co-dimension $r$
which satisfies the following two conditions:
 
1). for each point $t\in U$, the intersection ${\mathcal
M}\cap N
\times t$ is a connected, compact submanifold of $N\times t$ of dimension
$n$.
 
2). for each point $p\in{\mathcal M}$, there exist $r$ holomorphic
functions 
$f_1=f_1(w,t),$$\cdots,$$f_r=f_r(w,t)$ defined on a neighborhood
${\mathcal U}_p$of $p$ in $N\times U$ such that
\[
rank\,\frac{\pa(f_1,\cdots, f_r)}{\pa(w^1,\cdots, w^{r+d})}=r
\]
and in ${\mathcal U}_p$, the submanifold ${\mathcal M}$ is defined by the 
simultaneous equations
\[
f_1(w,t)=f_2(w,t)=\cdots =f_r(w,t)=0
\]
We call $U$ the parameter manifold or the base space of the
family $(N,U)$. We denote the family $(N,
U)$ simply by $N$ when we need not indicate 
the base space
$U$. For each point $t\in U$, we set
\[
M_t\times t={\mathcal M}\cap N\times t
\]
the submanifold $M_t$ of $N$ thus defined will be called 
the fiber of ${\mathcal M}$ over
$t$. We may identify $M_t\times t$ with
$M_t$ and consider $M_t$ as a 
{\it family consisting of compact submanifold} $M_t$,
$t\in U$ of $N$.
 \end{definition}

\begin{definition}
We say $({\frk X}, U)$ is the local total  family of $M$, if 
 $U$ is a neighborhood of $C^d$ with $d=\dim\, H^1(M,T_M)$ 
where $T_M$ is the holomorphic tangent bundle of $M$, and a
projection
\[
\pi: {\frk X}\rightarrow U
\]
such that it is holomorphic, surjective, of rank $d$ 
and such that for all $t\in U$,
$\pi^{-1}(t)$ is a deformation of complex structure of the center fiber
$\pi^{-1}(0)=M$.
\end{definition}

We  state the main result of this paper.

\begin{theorem}
Suppose $M$ is a compact complex manifold whose deformation of its
complex structure is unobstructed. Let $N$ be another complex manifold.
Suppose that
\[
i: M\rightarrow N
\]
is a holomorphic embedding. If $H^1(M, T_N|_M)=0$ where
$T_N|_M$ is the restriction of the holomorphic
tangent bundle of $N$ to $M$, then there is a
holomorphic
map
\[
f:{\frk X}\rightarrow N
\]
such that $f|_M=i$. 
Here $({\frk X},U)$ is the local total family of $M$.
Furthermore, $f|_{\pi^{-1}(t)}$ is an embedding of
$\pi^{-1}(t)$ to $N$ for $t\in U$.
\end{theorem}

{\bf Proof:} 
Suppose $\cup_{j\in I} U_j\supset {\frk X}$
is an open covering. On each $U_j$, $j\in I$, 
suppose $(z_j^1,\cdots ,z_j^n,t)$
is a local
coordinate such that
\[
\pi(z_j^1,\cdots ,z_j^n,t)=t,\qquad t\in U\subset C^d
\]
It is obvious that for fixed $t$, $(z_j^1,\cdots z_j^n)$ will be local
coordinate
for $\pi^{-1}(t)=M_t$. 
We further assume that $U_j$ is defined by
\[
U_j=\{ |z_j|=Max_\alpha |z_j^\alpha|<1\}
\]
We have, however, 
holomorphic functions $g_{jk}$ such that
\[
z_j^{\alpha}=g_{jk}^\alpha(z_k,t),\qquad j,k\in I
\]
for $\alpha=1,\cdots,n$ and for  $U_j\cap U_k\neq \emptyset$.

Now we suppose $\cup V_A\supset N$ is an open covering of $N$. Suppose 
\[
i(U_j)\subset V_{A(j)}, \qquad j\in I
\]
for some $A(j)$ of $j$.
Suppose $(w_A^1,\cdots,w_A^r)$
 is the local holomorphic coordinate  chart of $N$ on $V_A$.
And
we have transition functions $h_{AB}$ on $V_A\cap V_B\neq\emptyset$,
\[
w_A^s=h^s_{AB}(w_B)
\]
for $s=1,\cdots r+n=\dim N$. And again, we assume 
\[
V_j=\{ |w_j|=~Max |w_j^s|<1\}
\]

For the sake of simplicity, we denote
$j$ for $A(j)$. In order to construct $f$, we need only have to construct
holomorphic
mappings $(f_j)$ such that
\[
f_j: U_j\rightarrow V_j,\qquad j\in I
\]
satisfying
\begin{equation}\label{8}
h^s_{jk}(f_k(z_k,t))=f^s_j(g_{jk}(z_k,t),t)\qquad 
on \quad U_j\cap U_k\neq\emptyset
,\quad j,k\in I
\end{equation}
for $s=1,\cdots,r+n$.

We set up some notations. Let
\[
f_k(z,t)=f_{k|0}+f_{k|1}+\cdots+f_{k|m}+\cdots
\]
be the decomposition of $f_k(z_k,t)$ into homogeneous polynomials of $t$ of
degree $m$. 
Of course, each $f_k$'s and $f_{k|m}$'s are vector valued functions
$f_k=(f_k^s), f_{k|m}=(f_{k|m}^s), s=1,\cdots,r+n$.
Suppose
\[
f_k^m=f_{k|0}+\cdots + f_{k|m},\qquad k\in I
\]
and let $a\equiv_m b$ means $a-b$ is of polynomial of $t$ of degree bigger than
or equal to $m+1$.

We construct $f_{k|m}$ inductively. First, set
\[
f_{j|0}=i_j(z_j)\qquad j\in I
\]
It is easy to check that
\[
h_{jk}(f_k^0(z_k,t))\equiv_0 f_j^0(g_{jk}(z_k,t),t)
\qquad U_j\cap U_k\neq\emptyset
\]
where $h_{jk}=(h_{jk}^s)_{s=1,\cdots,r+n}$.

Now suppose for integer $m$, $f_k^m$ is constructed, and
\begin{equation}\label{7}
h_{jk}(f_k^m(z_k,t))\equiv_mf_j^m(g_{jk}(z_k,t),t)
\qquad U_j\cap U_k\neq\emptyset
\end{equation}
Define
\begin{equation}\label{1}
\Psi_{jk}(z_k,t)\equiv_{m+1}h_{jk}(f_k^m(z_k,t))-f_j^m(g_{jk}(z_k,t),t)
\quad U_j\cap U_k\neq\emptyset
\end{equation}

Then

{\bf Claim:}
\[
\Psi_{ik}=\Psi_{ij}+\frac{\pa w_i}{\pa w_j^s}|_{t=0}\Psi^s_{jk}
\qquad 
on\quad U_i\cap U_j\cap U_k\neq\emptyset
\]

{\bf Proof of the Claim:} By Equation~(\ref{1}), we have
\[
f_i^m(g_{ik}(z_k,t),t)\equiv_{m+1}h_{ik}(f_k^m(z_k,t))-\Psi_{ik}(z_k,t)
\qquad U_i\cap U_k\neq\emptyset
\]

Thus on $U_i\cap U_j\cap U_k\neq\emptyset$, we have
\begin{align*}
&h_{ij}(f_j^m(g_{jk}(z_k,t),t))\equiv_{m+1}
h_{ij}(h_{jk}(f_k^m(z_k,t))-\Psi_{jk}(z_k,t))\\
& \equiv_{m+1}h_{ij}(h_{jk}(f_k^m(z_k,t)))-\frac{\pa w_i}{\pa
w^s_j}\Psi^s_{jk}(z_k,t)\\
&\equiv_{m+1}h_{ik}(f_k^m(z_k,t))-\frac{\pa w_i}{\pa w_j^s}\Psi^s_{jk}(z_k,t)
\end{align*}

Note that $\Psi^s_{jk}(z_k,t)\equiv_m0$. So
\[
\frac{\pa w_i}{\pa w_j^s}\Psi_{jk}^s(z_k,t)\equiv_{m+1}
\frac{\pa w_i}{\pa w_j^s}|_{t=0}\Psi_{jk}^s(z_k,t)
\]

We have
\[
\Psi_{ij}(g_{jk}(z_k,t),t)\equiv_{m+1}\Psi_{ik}(z_k,t)
-\frac{\pa w_i}{\pa w_j^s}|_{t=0}\Psi^s_{jk}(z_k,t)
\qquad U_j\cap U_k\neq\emptyset
\]

On the other hand
\[
g_{jk}(z_k,t)\equiv_0 z_j
\qquad U_j\cap U_k\neq\emptyset
\]

So

\[
\Psi_{ij}(z_j,t)\equiv_{m+1}\Psi_{ij}(g_{jk}(z_k,t),t)
\]

Thus  we have
\begin{equation}\label{5}
\Psi_{ij}(z_j,t)=\Psi_{ik}(z_k,t)-\frac{\pa w_i}{\pa
w_j^s}|_{t=0}\Psi^s_{jk}(z_k,t)\qquad 
U_j\cap U_k\neq\emptyset
\end{equation}
The claim is proved.\qed

\smallskip

\smallskip

Suppose ${\mathcal U}=(U_j)$ is the covering of $M$, we see
from Equation~(\ref{5})
$\{\Psi_{ij}\}_{i,j\in I}$ defined a cocycle of $Z^1({\mathcal
U},T_N|_M)$. Thus $\{\Psi_{ij}\}_{i,j\in I}\in H^1({\mathcal U}, T_N|_M)$.
With a good covering ${\mathcal U}$, we have $H^1({\mathcal
U},T_N|_M)=H^1(M,T_N|_M)$ and by the assumption, the latter is zero.
So we can find $\{\Psi_j\}$ such that
\[
\Psi_{jk}=\Psi_j-\frac{\pa w_j}{\pa w_k^s}|_{t=0}\Psi^s_k\qquad 
U_j\cap U_k\neq\emptyset
\]

We then define $f_k^{m+1}$ inductively as
\begin{equation}\label{2}
f_k^{m+1}(z_k,t)=f_k^m(z_k,t)+\Psi_k(z_k,t),\qquad k\in I
\end{equation}

With this definition, we have
\begin{align*}
&h_{jk}(f_k^{m+1}(z_k,t))\equiv_{m+1}h_{jk}(f_k^m(z_k,t)+\Psi_k(z_k,t))\\
&\equiv_{m+1}h_{jk}(f_k^m(z_k,t))+\frac{\pa w_j}{\pa
w_k^s}\Psi^s_k(z_k,t)\\
&\equiv_{m+1}\Psi_{jk}(z_k,t)+\frac{\pa w_j}{\pa w^s_k}\Psi^s_k(z_k,t)
+f_j^m(g_{jk}(z_k,t),t)\\
&\equiv_{m+1}\Psi_j(z_j,t)+f_j^m(g_{jk}(z_k,t),t)\\
&\equiv_{m+1} f_j^{m+1}(g_{jk}(z_k,t),t)
\end{align*}

Now we have got a formal series
\[
f_{k|0}+f_{k|1}+\cdots +f_{k|m} +\cdots
\qquad k\in I
\]
which satisfies Equation~(\ref{7}) for any $m$.
If it converges, then $f_k$'s and  $f$ will be holomorphic
and $f$ will satisfy Equation~(\ref{8}). We put
off the proof of the convergence to the next section. At this moment, we
assume the convergence is true.

By continuity, for fixed $t$, $t$ being sufficiently small, $f$ will be an
immersion on $\pi^{-1}(t)$. We claim that for sufficiently small $t$, $f$
is an embedding. Suppose not, then we can find $t_n\rightarrow 0$ and
$x_n,y_n\in \pi^{-1}(t_n)$ such that $x_n\neq y_n$ but
$f(x_n,t)=f(y_n,t)$.
Suppose  $x_n\rightarrow x$ and $y_n\rightarrow y$. We see $x=y$,
otherwise it
will contradict to the fact that $i$ is an embedding. But if $x=y$, it
will contradict to the fact that $f$ is an immersion on each fiber.

Now we give an important example of manifolds such that
\linebreak
$H^1(M,T_N|_M)=0$.

\begin{prop}
Let $M$ be a simply connected Calabi-Yau threefold. $M\rightarrow CP^D$ is
an embedding. Then we have
\[
H^1(M,T_{CP^D}|_M)=0
\]
\end{prop}

{\bf Proof:}
From the Euler exact sequence
\begin{equation}\label{x}
0\rightarrow C\rightarrow \oplus_{D+1}{\mathcal O}(1)\rightarrow
T_{CP^D}\rightarrow 0
\end{equation}
We have the long exact sequence
\[
\cdots\rightarrow H^1(M,\oplus_{d+1}{\mathcal O}(1))
\rightarrow H^1(M,T_{CP^D}|_M)\rightarrow
H^2(M,{\mathcal O})\rightarrow\cdots
\]
By Kodaira Vanishing theorem
\[
H^1(M,\oplus_{D+1}{\mathcal O}(1))=0
\]
By Dolbeault theorem and Serre Duality $H^2(M,{\mathcal O})=0$, we have
\[
H^1(M,T_{CP^D}|_M)=0
\]

\begin{cor}
Suppose $M$ is a simply connected Calabi-Yau threefold. If $M$ is embedded
to some $CP^D$, then $M_t=\pi^{-1}(t)$ can also be embedded to the same
$CP^D$ for small $t$.
 \end{cor}

The following example showed that in general,
$H^1(M, T_{CP^D}|_M/T_M)\neq 0$.

{\bf Example.\,\,} Suppose ${\mathcal N}=T_{CP^D}|_M/T_M$ is the normal
bundle
of $M$ in
$CP^D$ in the previous proposition. Then in general, $H^1(M,{\mathcal
N})\neq 0$.

\smallskip

\smallskip

{\bf Proof:} We have the exact sequence:
\[
0\rightarrow T_M\rightarrow T_{CP^D}|_M
\rightarrow {\mathcal N}|_M\rightarrow 0
\]
So the long exact sequence gives
\begin{align*}
&\cdots\rightarrow
H^1(M,T_{CP^D}|_M)\rightarrow H^1(M,{\mathcal N})
\rightarrow H^2(M,T_M)\\
&\rightarrow H^2(M,T_{CP^D}|_M)\rightarrow \cdots
\end{align*}

However, by Serre Duality, $H^2(M,T_M)=H^{1,1}(M)$. And it is easy to see
from the Euler
Sequence~(\ref{x}) that $\dim H^2(M,T_{CP^D}|_M)=1$. Thus in general
$\dim H^1(M,{\mathcal N})\geq \dim H^{1,1}(M)-1$ and is not zero.

\section{The Convergence}
Now we shall show that the power series $\sum f_k(z_k,t)$
in Equation~(\ref{2})
 for all $k\in I$,
converge
for $|t|<\eps_0$, $\eps_0$ being a sufficiently small positive number,
provided
that we choose for each $m$ the homogeneous polynomials $f_{k|m+1}(z_k,t), k\in
I$, satisfying in a proper manner.

For any vector $\xi=(\xi^1,\xi^2,\cdots,\xi^\lambda,\cdots)$, we define
\[
|\xi|=Max_\lambda|\xi^\lambda|
\]

Consider a power series
\[
\xi(z,u)=\sum\xi_{lm,\cdots,n}u^l_1u^m_2\cdots u^n_q
\]
in $u_1,\cdots,u_q$ whose coefficients $\xi_{lm,\cdots,n}$ are vector 
valued
functions of $z$ and a power series
\[
a(u)=a_{lm,\cdots,n}u_1^lu_2^m\cdots u_q^n,\qquad a_{lm,\cdots,n}\geq 0
\]

We indicate by writing $\xi(z,u)\ll a(u)$ that
\[
|\xi_{lm,\cdots,n}(z)|\leq a_{lm,\cdots,n}
\]

Let
\[
A(t)=\frac{a}{16b}\sum_{n=1}^{\infty}\frac{1}{n^2}b^n(t_1+\cdots+t_l)^n
\]
where $a$ and $b$ are positive constants. We have
\[
A(t)^\gamma\ll (\frac ab)^{\gamma-1}A(t),\qquad \gamma=2,3,\cdots
\]

For our purpose it suffices to prove the inequalities
\[
f_k(z_k,t)\ll A(t),\qquad i\in I
\]

In what follows we denote by $c_0,c_1,c_2,\cdots$ positive constants {\it which
are greater than 1}. We may assume that
\[
|\frac{\pa h_{ij}^\lambda}{\pa w_j^{\mu}}|<c_0,\qquad  c_0>1
\]

For the sake of simplicity, we denote $f_k^m(z_k,t)-i(z_k)$ by
$f_k^m(z_k,t)$. Then for sufficiently large $a$, we have

\[
f_k^1(z_k,t)\ll \frac{a}{16}(t_1+\cdots+t_d)\ll A(t),\qquad k\in I
\]

Now assuming the inequalities
\[
f_k^m(z_k,t)\ll A(t),\qquad  k\in I
\]
for an integer $m\geq 1$, we shall estimate the coefficients of the
homogeneous
polynomials $\Psi_{ik}(z,t)$. We expand $h_{ik}$ and $g_{ik}$ into
power series,
whose coefficients are vector values holomorphic functions:

\begin{align*}
&h_{ik}(w_k)\ll\sum_{\alpha=0}^{\infty}c_1^\alpha(w_k^1+\cdots
w_k^{r+n})^\alpha\\
&g_{ik}(z_k)\ll\sum_{\alpha=0}^{\infty}c_1^\alpha(z_k^1+\cdots+z_k^n)^\alpha
\end{align*}
for $i,k\in I$.

Recall that in Equation~(\ref{1})
\[
\Psi_{jk}=[h_{jk}(f_k^m(z_k,t))]_{m+1}-[f_j^m(g_{jk}(z_k,t),t)]_{m+1}
\]
where $[a]_{m+1}$ is of polynomial of $t$ of degree bigger than $m$.
First we estimate $[h_{jk}(f_k^m(z_k,t))]_{m+1}$. The terms which are
linear
in $h_{jk}$ contributes nothing to $[h_{jk}(f_k^m(z_k,t))]_{m+1}$.
So we have
\begin{align*}
&[h_{jk}(f_k^m(z_k,t))]_{m+1}\ll\sum_{\alpha=2}^{\infty}
c_1^\alpha (r+n)^\alpha A(t)^\alpha\\&
\ll c_1(r+n)A(t)\sum_{\alpha=1}^{\infty}(\frac{c_1(r+n)a}{b})^\alpha
\end{align*}

Assuming that
\[
b>2c_1(r+n)a
\]
we obtain therefor
\begin{equation}\label{21}
[h_{jk}(f_k^m(z_k,t))]_{m+1}\ll 2c_1^2r^2ab^{-1}A(t)
\end{equation}

On the other hand
\[
[f_j^m(g_{jk}(z_k,t),t)]_{m+1}
=[f_j^m(g_{jk}(z_k,t),t)-f_j^m(z_j,t)]_{m+1}
\]
Denote by $U_i^\delta$ the subdomain of $U_i$ consisting of all points
$z_j=(z_j,0)$, $|z_j|<1-\delta$. We fix a positive number $\delta$ such
that
$\{U^\delta_i|i\in I\}$ forms a covering of $M$. Take a point $z\in
U_k\cap U_i^\delta$ and let $z_k$ and $z_j$ be the local coordinates of
$z$
on $U_k$ and $U_j$ respectively. Obviously, we have
\[
z_j=g_{jk}(z_k,0),\qquad |z_k|<1, |z_j|<1-\delta
\]

Letting $y=(y_1,\cdots, y_n)$, we expand the coefficients of polynomial
$f_i^m(z_j+y,t)$ into power series. Suppose $|y|<\delta$, we have
\[
[f_j^m(z_j+y,t)-f_j^m(z_j,t)]_{m+1}
\ll
A(t)(\Pi_{\alpha=1}^n(1-\frac{|y_\alpha|}{\delta})^{-1}-1)
\]

Now if $t$ is small, and 
$\mu=Max_{j,k} |g_{jk}(z_k,t)-z_j|<\delta$, then
\begin{equation}\label{3}
[f_j^m(g_{jk}(z_k,t),t)]_{m+1}\ll ( (1-\frac{\mu}{\delta})^{-n}-1)A(t)
\end{equation}

So from Equation~(\ref{21}) and ~(\ref{3}), we have

\[
\Psi_{jk}\ll (2c_1^2(r+n)^2ab^{-1}+((1-\frac{\mu}{\delta})^{-n}-1))A(t)
\]

We take arbitrary point $z\in U_k\cap U_i$ and choose a domain
$U_j^\delta$, which contains $z$, then 
\[
\Psi_{ik}=\frac{\pa w_i}{\pa w_j}\Psi_{jk}
-\frac{\pa w_i}{\pa w_j}\Psi_{ji}
\]
Thus
\[
\Psi_{ik}\ll 2c_0
(2c_1^2(r+n)^2ab^{-1}+((1-\frac{\mu}{\delta})^{-n}-1))A(t)
\]

Let
$c_3=2c_0(2c_1^2(r+n)^2ab^{-1}+((1-\frac{\mu}{\delta})^{-n}-1))$. 
Then
\begin{equation}\label{4}
\Psi_{jk}\ll c_3A(t)\qquad z\in U_k\cup U_i
\end{equation}

The following lemma can be proved by an elementary consideration.
~\cite{KS1}.

\begin{lemma}
We can choose the homogeneous polynomials $\Psi_i, i\in I$
satisfying
\[
\Psi_{ik}=\Psi_i-\frac{\pa w_i}{\pa w_k^s}\Psi_k^s
\]
in such a way that
\[
\Psi_i\ll c_3c_4A(t)
\]
where $c_4>1$ is a constant independent of $m$.
\end{lemma}\qed

Since $\mu$ is independent of $m$, and $\mu\rightarrow 0$ as $t\rightarrow
0$, we can choose $t$ small enough such that
$2c_0((1-\mu/\delta)^{-n}-1)<\frac{1}{2}$. Choosing $b$ large enough
so that
$4c_0c_1^2(r+n)^2ab^{-1}<\frac 12$, we have
\[
f^{m+1}_j\ll A(t),\qquad
j\in I
\]
The convergence of the power series is proved.

\end{document}